\documentclass[a4paper,11pt]{amsart}

\usepackage{amssymb,amsmath,eufrak,bbm,hyperref}
\usepackage{graphicx}
\usepackage{subfigure}
\usepackage{epsfig}
\usepackage{stmaryrd}
\author[Florent Benaych-Georges]{Florent Benaych-Georges}\address{Florent Benaych-Georges, LPMA,  UPMC Univ Paris 6, Case courier 188, 4 Place Jussieu, 75252 Paris Cedex 05, France and CMAP, \'Ecole Polytechnique, route de Saclay, 91128 Palaiseau Cedex, France.}\email{florent.benaych@upmc.fr}

\author[Nathana\"el Enriquez]{Nathana\"el Enriquez}\address{Nathana\"el Enriquez, MODAL'X, 200 Avenue de la R\'epublique, 92001, Nanterre and LPMA,   Case courier 188, 4 Place Jussieu, 75252 Paris Cedex 05, France.} \email{nenriquez@u-paris10.fr}
\title{Perturbations of  diagonal matrices by band random matrices}
\date{\today}

\thanks{This work was partially supported by the \emph{Agence Nationale de la Recherche} grant ANR-08-BLAN-0311-03.}

\newcommand{\ps}{\f{\partial}{\partial s}}
\newcommand{\pt}{\f{\partial}{\partial t}}
\newcommand{\one}{\mathbbm{1}}

\newcommand{\rfl}{\rfloor}
\newcommand{\lfl}{\lfloor}
\newcommand{\si}{\sigma}

\newcommand{\dens}{\operatorname{density}}
\newcommand{\splaw}{\operatorname{spectral \; law}}

\newcommand{\diag}{\operatorname{diag}}

\newcommand{\ninf}{\underset{n\to\infty}{\longrightarrow}}

\newcommand{\E}{\mathbb{E}}

\newcommand{\R}{\mathbb{R}}
\newcommand{\C}{\mathbb{C}}

\newcommand{\ud}{\mathrm{d}}

\newcommand{\pro}{probability }

\newcommand{\f}{\frac}
\newcommand{\ff}{\frac{1}}
\newcommand{\lf}{\left}
\newcommand{\ri}{\right}

\newcommand{\st}{such that }
\newcommand{\la}{\lambda}
\newcommand{\La}{\Lambda}

\newcommand{\ste}{\, ;\, }
\newcommand{\mc}{\mathcal }

\newcommand{\eps}{\varepsilon}

\newcommand{\bxp}{\boxplus}

\newcommand{\al}{\alpha}
\newcommand{\tta}{\theta}

\newcommand{\bbm}{\begin{bmatrix}}
\newcommand{\ebm}{\end{bmatrix}}
\newcommand{\bes}{\begin{equation*}}
\newcommand{\ees}{\end{equation*}}
\newcommand{\be}{\begin{equation}}
\newcommand{\ee}{\end{equation}}
\newcommand{\beqy}{\begin{eqnarray}}
\newcommand{\eeqy}{\end{eqnarray}}
\newcommand{\beq}{\begin{eqnarray*}}
\newcommand{\eeq}{\end{eqnarray*}}
\newcommand{\lto}{\longrightarrow}
\newcommand{\sd}{\,\cdot\,}

\newtheorem{Th}{Theorem}

\newtheorem{rmq}{Remark}

\newenvironment{pr}{\noindent {\bf Proof. }}{\hfill$\square$}

\long\def\symbolfootnote[#1]#2{\begingroup
\def\thefootnote{\fnsymbol{footnote}}\footnote[#1]{#2}\endgroup}

\begin{document}
\maketitle
\begin{abstract}We exhibit an explicit formula for the spectral density of a (large) random matrix which is a  diagonal matrix whose spectral density converges, perturbated by the addition of a symmetric matrix with Gaussian entries and a given (small) limiting variance profile.
\end{abstract}
\section{Perturbation of the spectral density of a large diagonal matrix}
In this paper, we consider the spectral measure of a   random matrix $D_n^\eps$ defined by $D_n^\eps=D_n+\sqrt{\f{\eps}{n}}X_n$, for $D_n$ a deterministic diagonal matrix whose spectral measure converges and $X_n$ an Hermitian or real symmetric matrix whose entries are  Gaussian independent variables, with a limiting variance profile (such matrices are called {\it band matrices}). We give a first order Taylor expansion, as $\eps\to 0$, of the limit spectral density, as $n\to\infty$,  of $D_n^\eps$.

The proof is elementary and based on a formula given in \cite{Sh96}   for the Cauchy transform of the  limit spectral distribution of $D_n^\eps$ as $n\to\infty$. 
\vskip1cm
 For each $n$, we consider an Hermitian or real symmetric random matrix $X_n=[x^n_{i,j}]_{i,j=1}^n$   and a real diagonal matrix $D_n=\diag(a_n(1),\ldots, a_n(n))$. We suppose that:\\
 \begin{enumerate}\item[(a)] the entries $x^n_{i,j}$ of $X_n$ are independent (up to symmetry), centered, Gaussian with  variance denoted  by $\si_n^2(i,j)$,\\
  \item[(b)] for a certain bounded function $\sigma$ defined on $[0,1]\times [0,1]$ and a certain bounded real function $f$ defined on $[0,1]$, we have, in the $L^\infty$ topology, 
  $$\si_n^2(\lfl nx\rfl,\lfl ny\rfl)\ninf \si^2(x,y)\qquad\textrm{ and }\qquad a_n(\lfl nx\rfl)\ninf f(x),$$
  \item[(c)] the set of discontinuities of the function $\si$ is closed and intersects a finite number of times any vertical line of the square $[0,1]^2$.\\\end{enumerate}

For $\eps\ge 0$,   let us define, for all $n$, $$D_n^\eps=D_n+\sqrt{\f{\eps}{n}}X_n.$$

It is  known,  from  Shlyakhtenko in \cite[Th. 4.3]{Sh96} (see also \cite{AZ05}, which also provides a fluctuation result), that as $n$ tends to infinity, the spectral distribution of $D_n^\eps$ tends to a limit $\mu_\eps$ with Cauchy transform $$C_\eps(z)=\int_{x=0}^1C_\eps(x,z)\ud x,$$ where $C_\eps(x,\sd)$ is defined by the fact that it is analytic, maps the upper half-plane $\C^+$ into the lower one $\C^-$, and satisfies the relation \be\label{eqC_epsilon}C_\eps(x,z)=\ff{z-f(x)-\eps\int_{y=0}^1 \si^2(x,y)C_\eps(y,z)\ud y}.\ee 

Our goal here is to understand $\mu_\eps -\mu$ for small values of $\eps$. Let us introduce the set $\mc{T}$ of test functions we shall use here. We define $$\mc{T}=\lf\{t  \longmapsto \ff{z-t}\ste z\in \C^+\ri\}.$$

 Let us now define the {\it Hilbert transform}, denoted by $H[u]$,  of a   function $u$:$$H[u](s):=\operatorname{p.v.}\int_{t\in\R}\f{u(t)}{s-t}\ud t= \int_{y\in \R}\f{u(s-y)-u(s)}{y}\ud y.$$

Before stating our main result, let us make some   assumptions on the functions  $\si$ and $f$: \begin{enumerate}\item[(d)] the push-forward $\mu$ of the uniform measure on $[0,1]$ by the function $f$ has a  density $\rho$ with respect to the Lebesgue measure on $\R$,\\
 \item[(e)] there exists a symmetric    function $\tau(\sd,\sd )$ \st for all $x,y$, $\si^2(x,y)=\tau(f(x),f(y))$,\\
 \item[(f)]  there exist $\eta_0>0, \al>0$ and $C<\infty$  \st   for almost all $s\in \R$, for all $t\in [s-\eta_0, s+\eta_0]$,\quad $|\tau(s,t)\rho(t)-\tau(s,s)\rho(s)|\le C|t-s|^\al.$
 \end{enumerate}
  
  Note that by hypothesis (f) and by the boundedness of the function $f$,  the function $$s\longmapsto \rho(s)H[\tau(s,\sd)\rho(\sd)](s) $$ is well defined and compactly supported.
   
\begin{Th}\label{main-th}Under the   hypotheses (a) to (f), as $\eps\to 0$, for all $g\in \mc{T}$,  $$\int g(s)\ud \mu_\eps(s)=\int g(s)\ud \mu(s)-\eps\int g'(s)  F(s) \ud s+o(\eps),$$ with $F(s):=-  \rho(s)H[\tau(s,\sd)\rho(\sd)](s) .$ 

As a consequence, if  the function $F(\sd) $ has bounded variations, then $$\mu_\eps=\mu+\eps   \ud F+o(\eps).$$
\end{Th}

\begin{rmq}Roughly speaking, this theorem states that $$\lim_{\eps\to 0}\lim_{n\to\infty} \f{\splaw(D_n^\eps)\;-\;\splaw(D_n)}{\eps}=\ud F.$$
It would be interesting to let $\eps$ and $n$ tend to $0$ and $\infty$ together, and to find out the adequate rate of convergence to get a deterministic limit or non degenerated fluctuations. We are working on this question.\end{rmq}

\begin{rmq}
This result provides an analogue, for our random matrix model, of the  following formula about real random variables (valid when $Y$ is centered and independent of $X$):
$$\dens_{X+\sqrt{\eps} Y}(s)=\dens_X(s)+\eps\f{\E[Y^2]}{2}\dens_X''(s)+o(\eps).$$
\end{rmq}

\begin{rmq}
In the case where $X_n$ is a GUE or GOE matrix, the limiting spectral distribution of $D_n^\eps$ as $n\to\infty$ is the free convolution of the limiting spectral distribution of $D_n$ with a semi-circle distribution. Several papers are devoted to the study of qualitative properties (like regularity) of the free convolution     (see \cite{Biane2,BVReg,BA,BNRFC,BBGA09}). Besides, it has recently been proved that   type-B free probability theory  allows to give Taylor expansions, for small values of $t$, of the moments of $\mu_t\bxp\nu_t$ for two time-depending \pro measures $\mu_t$ and $\nu_t$ (see \cite{BS_typeB,Fevrier-NicaJFA,Fevrier-Higher-order}). Our work differs from the ones mentioned above by the fact that we allow to perturb $D_n$ by any band matrix, but also by the fact that it is focused on the density and not on the moments, giving an explicit formula rather than qualitative properties.
\end{rmq}

\begin{pr}For all $z\in \C^+$, we have  \be\label{majC_epsilon}\lf|C_\eps(x,z)\ri|\le \ff{\Im z} .\ee

Indeed, for all $y,z$ \st $z\in \C^+$, $C_\eps(y,z)\in \C^-$. As a consequence, the imaginary part of the denominator of the right hand term of \eqref{eqC_epsilon} is larger than $\Im(z)$.

Hence by \eqref{eqC_epsilon} and \eqref{majC_epsilon}, as $\eps\to 0$,  
$C_\eps(x, z)\lto \ff{z-f(x)}$ uniformly in $x$.


From what precedes,  \beq C_\eps(x,z)-\ff{z-f(x)}&=&\f{\eps\int_{y=0}^1\si^2(x,y)C_\eps(y,z)\ud y}{(z-f(x)-\eps\int_{y=0}^1\si^2(x,y)C_\eps(y,z)\ud y)(z-f(x))}\\ 
&=&
\eps\ff{(z-f(x))^2} \int_{y=0}^1\si^2(x,y)C_\eps(y,z)\ud y+o(\eps)\\
&=&
\eps\ff{(z-f(x))^2} \int_{y=0}^1\f{\si^2(x,y)}{z-f(y)}\ud y+o(\eps)
\eeq
where each $o(\eps)$   is   uniform in $x\in [0,1]$. 

But for all $a\ne b$, $\ff{(z-a)^2(z-b)}=\ff{(a-b)^2}\lf(\ff{z-b}-\ff{z-a}-\f{b-a}{(z-a)^2}\ri)$, hence since the Lebesgue measure of the set $\{y\in [0,1]\ste f(y)=f(x)\}$ is null, we have  \beq\ff{(z-f(x))^2} \int_{y=0}^1\f{\si^2(x,y)}{z-f(y)}\ud y&=& \int_{y=0}^1 \f{\si^2(x,y)}{(f(x)-f(y))^2}\lf(\ff{z-f(y)}-\ff{z-f(x)}-\f{f(y)-f(x)}{(z-f(x))^2} \ri)  \ud y.\eeq

As a consequence, it follows by an integration in $x\in  [0,1]$ that $$C_\eps(z)-C(z)=\eps\int_{x=0}^1\int_{y=0}^1 \f{\si^2(x,y)}{(f(x)-f(y))^2}\lf(\ff{z-f(y)}-\ff{z-f(x)}-\f{f(y)-f(x)}{(z-f(x))^2} \ri)  \ud y\ud x+o(\eps),$$where $C(\cdot)$ is the Cauchy transform of $\mu$. 

Let us now recall that the push-forward of the uniform law on $[0,1]$ by $f$ is the measure $\rho(x)\ud x$ and that $\si^2(x,y)$ can be rewritten $\si^2(x,y)=\tau(f(x),f(y))$. Hence  $$C_\eps(z)-C(z)=\eps\int_{s\in \R}\int_{t\in \R}\{\ff{z-t}-\ff{z-s}-\ff{(z-s)^2}(t-s)\}\f{\tau(s,t)}{(s-t)^2}\rho(s)\rho(t)\ud t\ud s+o(\eps).$$

This allows us to write that for any test function $g\in \mc{T}$,    $$\lim_{\eps\to 0} \f{\mu_\eps(g)-\mu(g)}{\eps}=\La(g),$$ where $$\La(g)=\int_{(s,t)\in \R^2} \{g(t)-g(s)-g'(s)(t-s)\}\f{ \tau(s,t)}{(t-s)^2}\rho(s)\rho(t)\ud t\ud s.$$ 
 Note that by the Taylor-Lagrange formula, for all $s,t$,$$ \lf| \{g(t)-g(s)-g'(s)(t-s)\}\f{ \tau(s,t)}{(t-s)^2}\rho(s)\rho(t) \ri|\le \f{ \rho(s)\rho(t) \times\|\tau(\sd,\sd)\|_{L^\infty} \|g''\|_{L^\infty}}{2},$$ so that, since $\rho$ is a density,   by     dominated convergence,   $$\La(g)=\lim_{\eta\to 0}\int_{\substack{(s,t)\in \R^2\\ |s-t|>\eta}}\{g(t)-g(s)-g'(s)(t-s)\}\f{ \tau(s,t)}{(t-s)^2}\rho(s)\rho(t)\ud s\ud t.$$ 
But by symmetry, 
 for all $\eta>0$, $$\int_{\substack{(s,t)\in \R^2\\ |s-t|>\eta}}\{g(t)-g(s)\}\f{ \tau(s,t)}{(t-s)^2}\rho(s)\rho(t)\ud s\ud t=0.$$As a consequence, 
 $\La(g)=\lim_{\eta\to 0}\La_\eta(g)$, with
$$\La_\eta(g):= \int_{\substack{(s,t)\in \R^2\\ |s-t|>\eta}}g'(s)\f{ \tau(s,t)}{s-t}\rho(s)\rho(t)\ud s\ud t.$$
 Let us prove that almost all $s\in \R$,  $\lim_{\eta\to 0}
 \int_{\substack{t\in \R\\ |s-t|>\eta}}\f{ \tau(s,t)\rho(s)\rho(t)}{s-t}\ud t$ exists and that $$\La(g)=\int_{s\in \R}g'(s) \lf(\lim_{\eta\to 0}
 \int_{\substack{t\in \R\\ |s-t|>\eta}}\f{ \tau(s,t)\rho(s)\rho(t)}{s-t}\ud t\ri)\ud s.$$
  For $\eta>0$ and $s\in \R$, set $$\tta_\eta(s):=\int_{\substack{t\in \R\\ |s-t|>\eta}}\f{ \tau(s,t)\rho(s)\rho(t)}{s-t}\ud t.$$Set  also $M:=\|f\|_{L^\infty}$. Then the  support of the function $\rho$ is contained in $[-M,M]$, and  so does  the support of the function  $\tta_\eta$, for any $\eta>0$.
  For almost all $s\in [-M,M]$, $\lim_{\eta\to 0}\tta_\eta(s)$ exists by the formula 
  $$\tta_\eta(s)=\int_{t\in [s-2M,s-\eta]\cup[s+\eta,s+2M]}\f{\tau(s,t)\rho(s)\rho(t)- \tau(s,s)\rho(s)\rho(s)}{s-t}\ud t$$ and by Hypothesis (f).
%
Moreover, for $\eta_0$ as in Hypothesis (f),
\beq \lf|\tta_\eta(s)\ri|
& \le& 2C\rho(s)\int_{t=s+\eta}^{s+\eta_0}(s-t)^{\al-1}\ud t+ \int_{t\in [s-2M,s-\eta_0]\cup[s+\eta_0,s+2M]}\f{  \tau(s,t)\rho(s)\rho(t)}{s-t}\ud t\\
&\le& \f{2C\rho(s)}{\al}(\eta_0)^\al+\ff{\eta_0}\int_{t\in \R} \tau(s,t)\rho(s)\rho(t)\ud s\ud t\\
&\le & \f{2C\rho(s)}{\al}(\eta_0)^\al+\f{\|\tau(\sd,\sd)\|_{L^\infty}}{\eta_0}\rho(s).
\eeq
Hence by dominated convergence, $\int_{s\in \R}g'(s)\lim_{\eta\to 0}\tta_\eta(s)\ud s=\lim_{\eta\to 0}\int_{s\in \R}g'(s)\tta_\eta(s)\ud s$, i.e.  $$\La(g)=\int_{s\in \R}g'(s) \lf(\lim_{\eta\to 0}
 \int_{\substack{t\in \R\\ |s-t|>\eta}}\f{ \tau(s,t)\rho(s)\rho(t)}{s-t}\ud t\ri)\ud s.$$
 
%
%
%
\end{pr}

\section{Examples}
\subsection{Perturbation of a uniform distribution by a standard band matrix}
\begin{figure}[h!]
\centering
\subfigure[Case where $n=4.10^3$, $\eps=10^{-2}$, with width $\ell=0.2$]{
\includegraphics[width=2.4in]{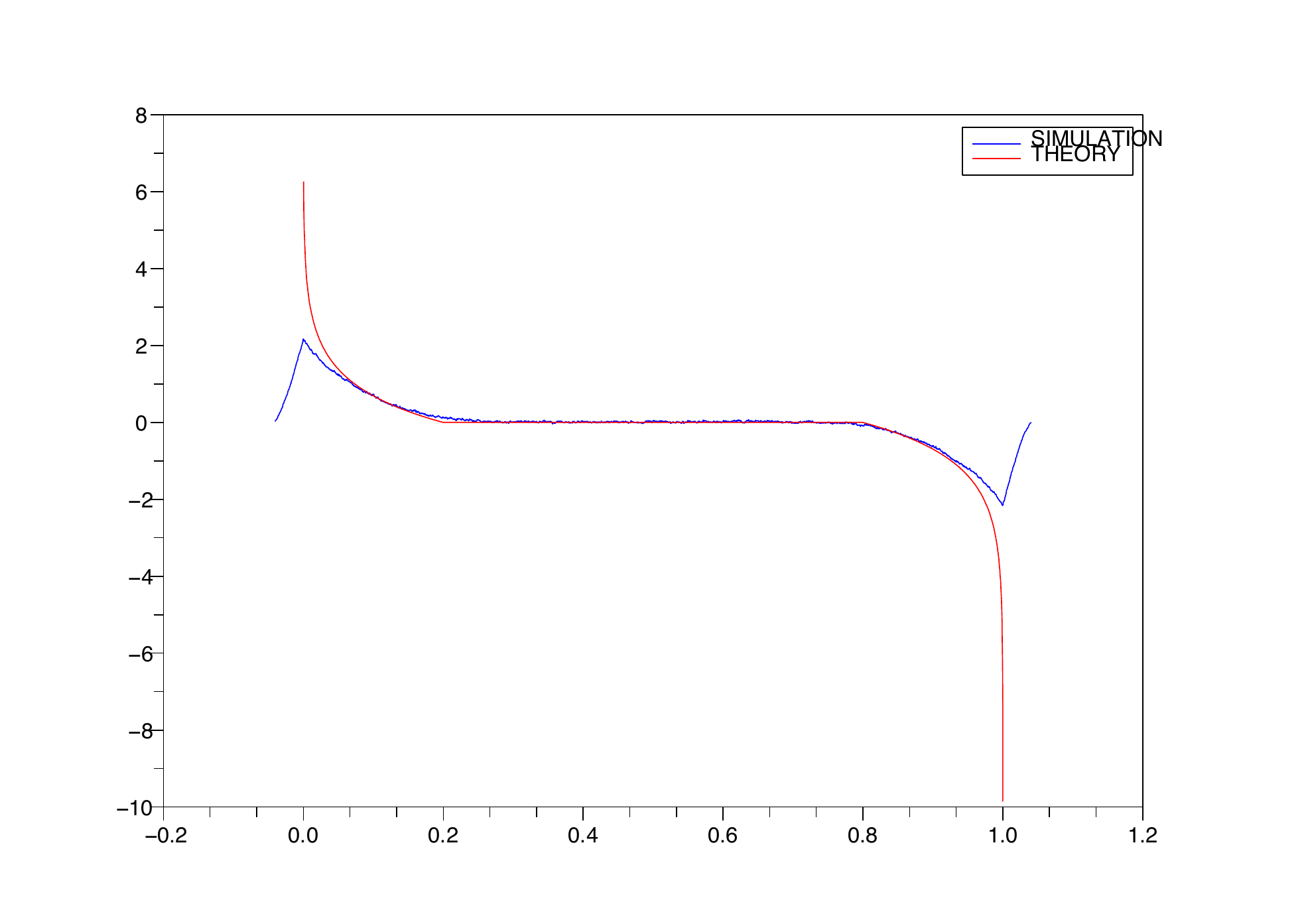}
\label{unif_width0.2}}\hspace{1.2in}
\subfigure[Case where $n=4.10^3$, $\eps=10^{-2}$, with width $\ell=0.9$]{
\includegraphics[width=2.4in]{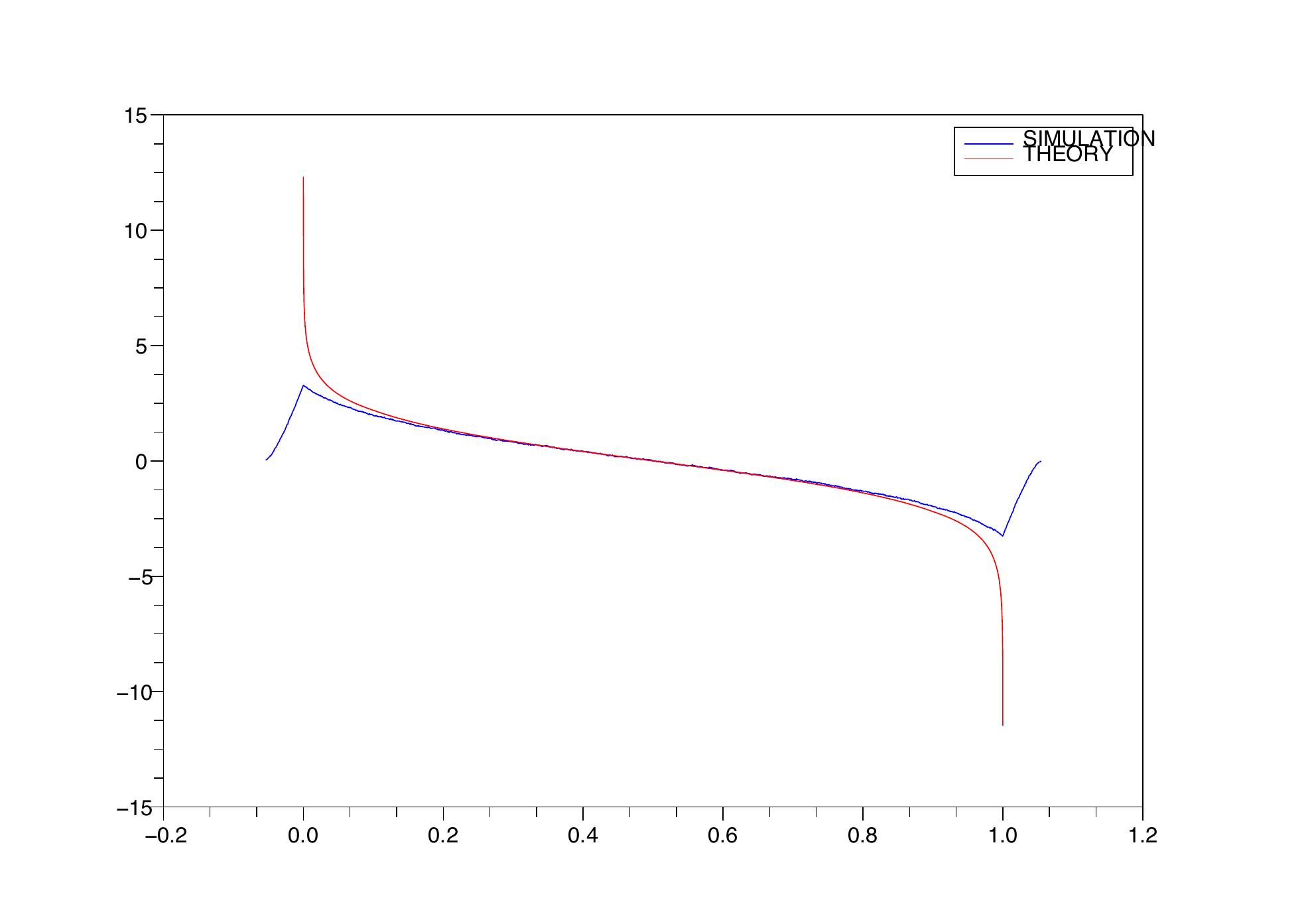}
\label{unif_width0.9}}
\caption{{\bf Perturbation of a uniform distribution by a standard band matrix:} plot   of the functions $F(\sd)$ and $\f{F_{D_n^\eps}(\sd)-F_{D_n}(\sd)}{\eps}$ (with $F_{D_n^\eps}(\sd)$ and $F_{D_n}(\sd)$ the cumulative eigenvalue distribution functions of $D_n^\eps$ and $D_n$) for different values of $\ell$.}
\label{fig:unif_band}
\end{figure}

Let us consider the case where $f(x)=x$ (so that $\mu$ is the uniform distribution on $[0,1]$) and $\si^2(x,y)=\one_{|y-x|\le \ell}$, where $\ell$ is a fixed parameter in $[0,1]$ (the width of the band).
In this case, $\tau(\sd,\sd)=\si^2(\sd,\sd)$ and $$F(s)=\one_{(0,1)}(s)\log\lf(\f{\ell \wedge(1-s)}{\ell\wedge s}\ri).$$
For small values of $\eps$ and large values of $n$, the density $\rho_\eps$ of the   eigenvalue  distribution $\mu_\eps$ of $D^\eps_n$ is approximately $$\rho_\eps(s)=\rho(s)+\eps\ps F(s)+o(\eps)=\one_{(0,1)}(s)-\eps \lf(\f{\one_{(0,\ell)}(s)}{s}+\f{\one_{(1-\ell,1)}(s)}{1-s}\ri)+o(\eps),$$ which means that the additive perturbation $\sqrt{\f{\eps}{n}}X_n$ alters the spectrum 
of $D_n$ essentially by decreasing the amount of extreme eigenvalues. This phenomenon is illustrated by Figure \ref{fig:unif_band} (where we ploted the cumulative distribution functions rather than the densities for visual reasons).

\subsection{Perturbation of the triangular pulse distribution by a GOE matrix}
\begin{figure}[h!]
\begin{center}
\includegraphics[width=14cm,height=4cm]{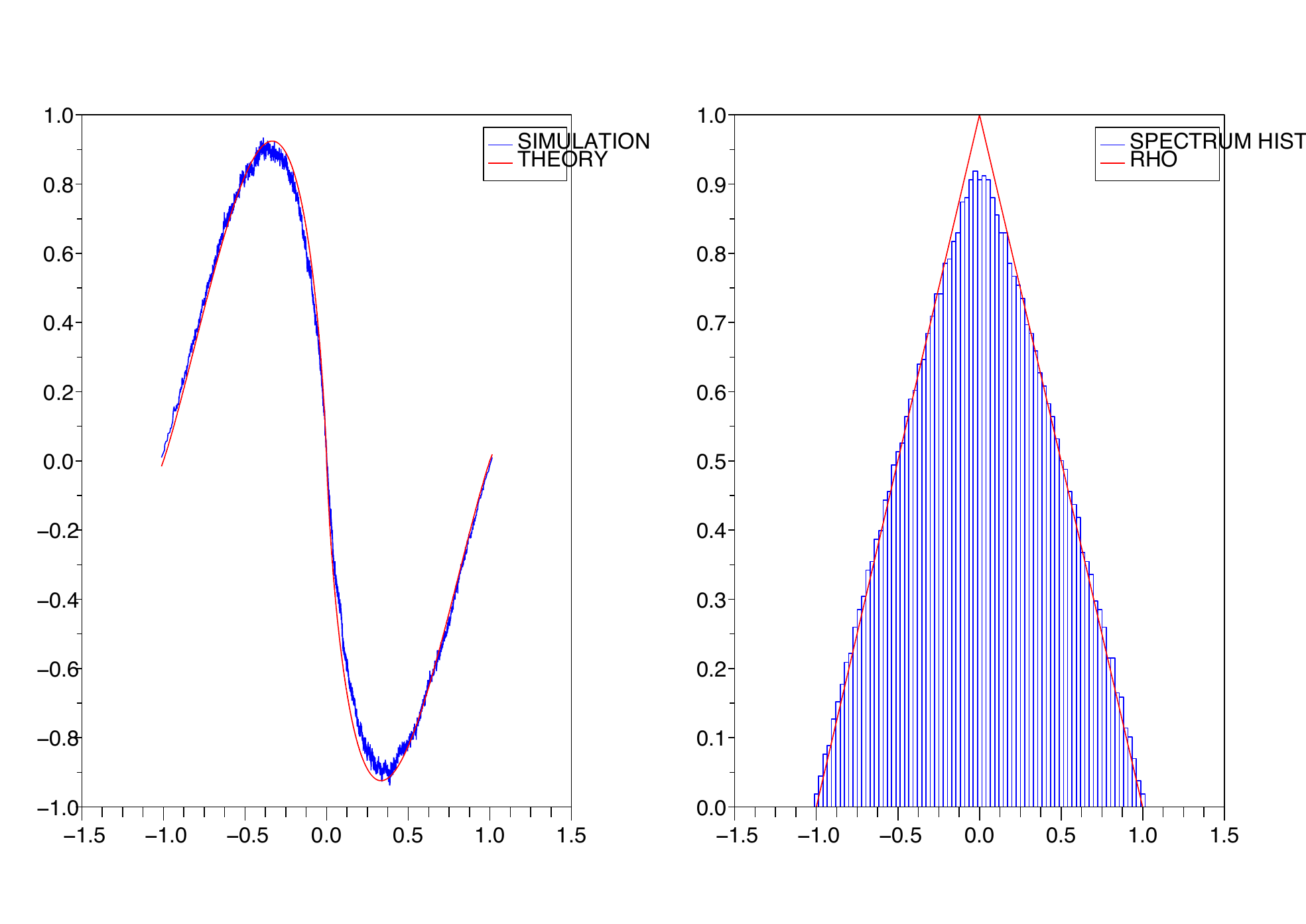}
\caption{{\bf Perturbation of the triangular pulse distribution by a GOE matrix:} {\it Left:} plot   of the functions $F(\sd)$ and $\f{F_{D_n^\eps}(\sd)-F_{D_n}(\sd)}{\eps}$ (with $F_{D_n^\eps}(\sd)$ and $F_{D_n}(\sd)$ the cumulative eigenvalue distribution functions of $D_n^\eps$ and $D_n$). {\it Right:} plot of the eigenvalues histogram of $D_n^\eps$ and of the spectral  density $\rho$ of $D_n$. On the right figure,  the (infinitesimal) increase of eigenvalues with respect to $\rho$  on $[-1,-0.5]\cup[0.5,1]$  and the (infinitesimal)  decrease  around zero can be observed, in agreement  with the fact that, as the left figure shows,   $F'\gg 0$ on (approximately) $[-1,-0.5]\cup[0.5,1]$ and $F'\ll 0$ around zero. Both figures were  made with the same simulation ($n=6.10^3$ and  $\eps=10^{-2}$).}\label{fig:tri_pluse}
\end{center}
\end{figure}

Let us consider the case where $\rho(x)=(1-|x|)\one_{[-1,1]}(x)$ and $\si^2\equiv 1$ (what follows can be adapted to the case $\si^2(x,y)=\one_{|y-x|\le \ell}$, but the formulas are a bit heavy).
In this case, thanks to the formula (9.6)  of $H[\rho(\sd)]$ given p. 509 of \cite{kingvol2}, we get $$F(s)=(1-|s|)\one_{[-1,1]}(s)\lf\{(1-s)\log(1-s)-(1+s)\log(1+s)+2s\log |s|\ri\}.$$
For small values of $\eps$ and large values of $n$, the density $\rho_\eps$ of the   eigenvalue  distribution $\mu_\eps$ of $D^\eps_n$ is approximately $$\rho_\eps(s)=\rho(s)+\eps\ps F(s)+o(\eps),$$  which implies that the additive perturbation $\sqrt{\f{\eps}{n}}X_n$ alters the spectrum 
of $D_n$   by increasing the amount of   eigenvalues in $[-1,-0.5]\cup[0.5,1]$ and decreasing the amount of eigenvalues around zero. This phenomenon is illustrated by Figure \ref{fig:tri_pluse}.

\subsection{Free convolution with a semi-circular distribution and complex Burger's equation}
Let us  consider the case where $\si^2\equiv 1$, which happens for example if the matrix $X_n$ is taken in the Gaussian Orthogonal Ensemble. In this case, by the theory of free probability developped by Dan Voiculescu (see e.g. \cite{Vo0} or \cite[Cor 5.4.11 (ii)]{alice-greg-ofer}), for all $t\ge 0$,  $$\mu_{{t}}=\mu\bxp\la_{t},$$
where $\la_t$ is the {\it semi-circular distribution with variance $t$}, i.e. the distribution with support $[-2\sqrt{t},2\sqrt{t}]$ and  density $\ff{2\pi t}\sqrt{4t-x^2}.$ In this case, we know by the work of Biane \cite[Cor. 2]{Biane2} that for all $t>0$, $\mu_t$ admits a density  $\rho_t$. By the  implicit function theorem, and the formula given in  \cite[Cor. 2]{Biane2}, one easily sees that the function $(s,t)\longmapsto \rho_t(s)$ is regular. 
Then, by Theorem \ref{main-th} and the fact that  the linear span of $\mc{T}$ is dense in the set of continuous functions on the real line with null limit at infinity, one easily recovers the following PDE, which is a kind of projection on the real axis of the imaginary part of   complex Burger's equation given in   \cite[Intro.]{Biane2}
\be\label{burgersreel}\begin{cases}\pt \rho_t(s)+\ps\{\rho_t(s)H[\rho_t(\sd)](s)\}=0,\\
\rho_0(s)=\rho(s).\end{cases}\ee

For example, 
if $\mu=\la_c$ for a certain $c>0$, then by the semi-group property of the semi-circle distribution \cite[Ex. 5.3.26]{alice-greg-ofer}, for all $t\ge 0$, $\mu_t=\la_{c+t}$ and $\rho_t(s)=\ff{2\pi(c+t)}\sqrt{4(c+t)-s^2}$. One can then verify \eqref{burgersreel}, using the formula (9.21) of $H[\rho_t(\sd)]$ given p. 511 of \cite{kingvol2}. 
\\

\noindent{\bf Acknowledgements.} It is a pleasure to thank Guy David for his useful advices about the Hilbert transform.

\end{document}